\documentclass{amsproc}
\usepackage{epsf}

\newcommand{\AB}{\A^\bullet}
\newcommand{\ABB}{\A^{\bullet,\bullet}}
\newcommand{\aut}{\operatorname{Aut}(G)}
\newcommand{\coh}{cohomology}
\newcommand{\Mc}{\overline{\mathcal{M}}}
\newcommand{\mgnb}[1]{\overline{\mathcal{M}}_{#1,n}}
\newcommand{\mgn}[1]{\mathcal{M}_{#1,n}}

\newcommand{\A}{\mathcal{A}}
\newcommand{\C}{\mathcal{C}}
\newcommand{\nc}{\mathbb{C}}
\newcommand{\gtm}{\mathfrak{M}}

\newtheorem{thm}{Theorem}[section]
\newtheorem{prop}[thm]{Proposition}
\newtheorem*{cor}{Corollary}

\theoremstyle{definition}

\theoremstyle{remark}

\newtheorem*{ack}{Acknowledgment}

\numberwithin{equation}{section}

\newcommand{\abs}[1]{\lvert#1\rvert}

\begin{document}

\title{Stability of the Rational Homotopy Type of Moduli Spaces}
\author{Alexander A. Voronov}
\address{Research Institute for Mathematical Sciences\\
Kyoto University\\
Sakyo-ku, Kyoto, 606-01\\Japan and Department of Mathematics\\ M.I.T.,
2-270\\ 77 Massachusetts Ave.\\ Cambridge, MA 02139-4307\\USA}
\email{voronov@kurims.kyoto-u.ac.jp
}
\thanks{Research supported in part by an AMS Centennial Fellowship.}

\date{August 20, 1997}

\dedicatory{Dedicated to Jim Stasheff on the occasion of his sixtieth
birthday.}

\begin{abstract}
We show that for $g \ge 2k+3$ the $k$-rational homotopy type of the
moduli space $\mgn{g}$ of algebraic curves of genus $g$ with $n$
punctures is independent of $g$, and the space $\mgn{g}$ is
$k$-formal. This implies the existence of a limiting rational homotopy
type of $\mgn{g}$ as $g \to \infty$ and the formality of it.
\end{abstract}

\maketitle

\section*{Introduction}

The description of the algebraic topology of the moduli space
$\mgn{g}$ of compact complex algebraic curves has long been a
tantalizing problem.  The idea of ``stable \coh '' of $\mgn{g}$ as the
genus $g \to \infty$, brought in by J.~L. Harer and D.~Mumford,
suggested a more graspable object to study as a first step. Mumford's
Conjecture \cite{mumford}, stating that the stable \coh\ is the
polynomial algebra on the so-called Mumford-Morita-Miller classes, has
become a new tantalizing problem since then. Recently M.~Pikaart
\cite{pikaart} has proven that the Hodge structure on the stable \coh\
is pure, thus providing more evidence to the conjecture, which easily
implies purity.

Our purpose in this note is to show that the rational homotopy type of
the moduli space $\mgn{g}$ stabilizes as $g \to \infty$ and prove
that the moduli space is formal in the stable limit. If the stable
\coh\ of the moduli space were just known to be the \coh\ of a certain
space $\mgn{\infty}$, then a theorem of E.~Miller and S.~Morita
\cite{miller} asserting that $H^\bullet (\mgn{\infty})$ is a free
graded commutative algebra would immediately imply the formality of
$\mgn{\infty}$. Thus, the essential result of this paper is showing
that the limit of the rational homotopy types of $\mgn{g}$ as $g \to
\infty$ exists. After this is done, the formality of the limit is
automatic.

Here are few words about the structure of the paper. Sections 1 and 2
are dedicated to the description of the (nonstable) rational homotopy
type of the moduli space $\mgn{g}$. Section~\ref{infty}, where the
stable rational homotopy type of $\mgn{g}$ is studied, is independent
of the first two sections.

It may be worth noting how this paper relates to Jim Stasheff's work.
Stasheff polyhedra $K_n$ are known for many remarkable properties.
One of them, that $K_n$ is a connected component of the real
compactification of the moduli space of real projective lines with
$n+1$ punctures, was pointed out by M.~Kontsevich in \cite{maxim}.
Thus, Stasheff polyhedra are real analogues of the compactified moduli
spaces $\mgnb{g}$. Or, perhaps, the moduli spaces $\mgnb{g}$ are
complex analogues of Stasheff polyhedra, the chicken and the egg
problem. The topology of $K_n$ is trivial ($K_n$ is contractible), and
the topology of $\mgn{g}$ and $\mgnb{g}$ is just the opposite. It is
rather the combinatorics of Stasheff polyhedra which makes them very
useful in the topology of loop spaces \cite{jim}. The combinatorics of
the genus zero spaces $\mgn{0}$ has similarly proven to be useful in
studying the topology of double loop spaces \cite{gj} and the
algebraic structure of 2d quantum field theory \cite{ksv1,ksv2,kvz}.
But what is the topology of these complex analogues of Stasheff
polyhedra, especially for a high genus?

We also implicitly use a result of S.~Halperin and Stasheff
\cite{halperin-stasheff} throughout the paper: formality over any
field of characteristic zero implies that over the rationals.

\begin{ack}
I am very grateful to Pierre Deligne, Dick Hain, Eduard Looijenga, and
Martin Pikaart for helpful discussions. I would like to thank the
Institute for Advanced Study in Princeton and the Research Institute
for Mathematical Sciences in Kyoto for their hospitality during the
summer 1997, when the essential part of the work on this paper was
done.
\end{ack}

\section{The rational homotopy type of a complex smooth
quasi-projective variety}

Here we recall J.~W. Morgan's description \cite{morgan} of a model of
the rational homotopy type of the complement $U$ of a normal crossings
divisor $D= \bigcup_{i=1}^r D_i$ in a compact complex manifold $X$.
The rational homotopy type will be understood in the sense of
D.~Sullivan \cite{sullivan}, see also P.~A. Griffiths and Morgan
\cite{griffiths-morgan} and D.~Lehmann \cite{lehmann}. It determines
the usual rational homotopy type of a topological space if it is
simply connected. In general it determines the rational nilpotent
completion thereof.

The rational homotopy type of $U$ is determined by the differential
graded (DG) commutative algebra that is nothing but the first term
$\AB = E_1^\bullet$ of the spectral sequence associated with the
weight filtration on the log-complex $\Omega^{\bullet,\bullet}_X(\log
D)$ of the smooth $(p,q)$-forms on $U$ with logarithmic singularities
along $D$. Otherwise, one can describe the same DG algebra as the
second term of the Leray spectral sequence of inclusion $U
\hookrightarrow X$. In any case, the model of $U$ is given by the DG
algebra $\AB$, where
\begin{align*}
\A^k & = \bigoplus_{p+q=k} \A^{p,q}, & \A^{p,q} & = \bigoplus_{
\substack{S \subset \{1,2,\dots,r\}\\
\abs{S} = -p}
}
H^{2p+q}(D_S, \nc),\\
p & \le 0, \; q \ge 0, & D_S & = \bigcap_{i \in S} D_i.
\end{align*}
The multiplication structure is given by
\[
a \cdot a' =
\begin{cases}
(-1)^{pq'+\epsilon} (a|_{D_S \cap D_{S'}}) \cup (a'|_{D_S \cap D_{S'}}) & \text{if $S \cap S' = \emptyset$},\\
0 & \text{otherwise}
\end{cases}
\]
for $a \in H^{2p+q}(D_S,\nc)$ and $a' \in H^{2p'+q'}(D_{S'}, \nc)$,
where $\epsilon$ is the sign of the shuffle putting the set $S \cap
S'$ in an increasing order, assuming each of the subsets $S$ and $S'$
to be already in an increasing order. Note that this multiplication
law makes $\ABB$ into a bigraded commutative algebra.

Finally, the differential $d: \A^{p,q} \to \A^{p+1,q}$ can be described as
\[
da = \sum_{j=1}^{-p} (-1)^{j-1} (\iota_j)_* a,
\]
where $a \in H^{2p+q}(D_S, \nc)$, $S = \{i_1, \dots, i_{-p}\}$, $\iota_j$ is
the natural embedding $D_S \subset D_{S\setminus i_j}$, and
$(\iota_j)_*: H^{2p+q}(D_S, \nc) \to H^{2p+q+2}(D_{S\setminus i_j},\nc)$ is the
Gysin map.

P.~Deligne \cite{del:hodge} proved that the spectral sequence
$E_1^{p,q} = \A^{p,q}$ degenerates at $E_2$, that is, the \coh\ of the
DG algebra $\AB$ is equal to $H^\bullet(U,\nc)$. The very DG algebra
$\AB$ describes the rational homotopy type of $U$, according to
Morgan's theorem \cite{morgan}.

\section{The rational homotopy type of $\mgn{g}$}
\label{mgn}

Results of the previous section apply to the moduli space $\mgn{g}$,
which is the complement of a normal crossing divisor in the
Deligne-Knudsen-Mumford compactification $\mgnb{g}$. The problem that
the space in question is a stack rather than variety does not arise,
because we work with the \textbf{complex coefficients}, as we will
assume throughout the paper. From the combinatorics of
Deligne-Knudsen-Mumford's construction, we can say more specifically,
cf.\ \cite{ksv2}, that
\[
\A^{p,q} = \bigoplus_G \; \left(
\bigotimes_{v \in G} H^{\bullet} (\Mc_{g(v), n(v)}) \right)_{2p+q}^{\aut},
\]
the summation running over all stable labeled $n$-graphs $G$ of genus
$g(G) = g$ and $v(G) = - p +1$ vertices. Here we refer to \emph{graphs}
of the following kind. Each graph is connected and has \emph{$n$
enumerated exterior edges}, edges which are incident with only one
vertex of the graph. Each vertex $v$ of the graph is \emph{labeled} by
a nonnegative integer $g(v)$, called the genus of a vertex. The
\emph{stability} condition means that any vertex $v$ labeled by
$g(v)=1$ should be incident with at least one edge (i.e., be at least
of valence one) and each vertex $v$ with $g(v) = 0$ should be at least
of valence three.  The \emph{genus} $g(G)$ of a graph $G$ is given by
the formula $g(G) = b_1(G) + \sum_v g(v)$, where $b_1(G)$ is the first
Betti number of the graph. The \emph{number $n(v)$} is the valence of
the vertex $v$ and $\aut$ is the \emph{automorphism group of a graph
$G$} (bijections on vertices and edges, preserving the exterior edges,
the labels of vertices and the incidence relation). The subscript
$2p+q$ in the formula refers to taking the homogeneous component of
this degree. The differential $d: \A^{p,q} \to \A^{p+1,q}$ is induced
by contracting interior edges in $G$, which corresponds to replacing a
neighborhood of a double point on a curve by a cylinder.

Our goal here is to look at the stable (as $g \to \infty$) \coh\ and
rational homotopy type of the moduli space $\mgn{g}$. The following
result of Pikaart \cite{pikaart} gives a certain clue on what is going
on within a ``stable range''.

\begin{thm}[Pikaart]
\label{pikaart}
The restriction mapping $H^k(\mgnb{g}) \linebreak[1]
\to
\linebreak[0]
H^k(\mgn{g})$ is surjective for $k \le (g-1)/2$.
\end{thm}

\begin{figure}[tb]
\centerline{\epsfxsize=1.5in \epsfbox{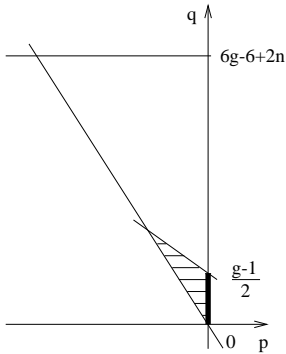}}
\caption{The algebra $\A^{p,q}$. The shaded region is the ``stable
range'', where the \coh\ is concentrated along the fat line.}
\label{graph}
\end{figure}

\begin{cor}
For $p+q \le (g-1)/2$, the \coh\ of $\A^{p,q}$ is nonzero only for
$p=0$, see Figure~$\ref{graph}$.
\end{cor}

\begin{proof}
By our construction of the spectral sequence, the natural composite
mapping $\A^{0,k} \to H^{0,k} (\ABB, d) \hookrightarrow H^k (\ABB, d)$
is the same as the restriction mapping $H^k(\mgnb{g}) \to
H^k(\mgn{g})$. Therefore, the $p \ne 0$ columns of $\A^{p,q}$ do not
contribute to the \coh\ of the DG algebra $\A$ in the stable range.
\end{proof}

\section{The stable limit}
\label{infty}

The stable limit in \coh\ of the moduli spaces $\mgn{g}$ is achieved,
roughly speaking, by gluing more and more handles to the complex
curve. This yields an inductive system (of isomorphisms) on the level
of \coh. But since there is no natural mapping between the moduli
spaces of different genera, one cannot speak of a limiting rational
homotopy type.  The question of taking the limit of the DG algebras
$\A$ of Section~\ref{mgn} for $\mgn{g}$ may not be so obviously
resolved, either, because these algebras are constructed out of the
\coh\ of $\mgnb{g'}$, $g'$ running between 0 and $g$. When $g\to
\infty$, $g'$ does not, and on top of that, taking the stable limit of
\coh\ of $\mgnb{g}$ requires a finer tuning, cf.\ Pikaart
\cite{pikaart}.

Our plan here is to show that a $k$-minimal model of $\mgn{g}$ is
independent of $k$ as long as $g \ge 2k+3$. In particular, the limit
of $k$-minimal models exists and may be called a ``$k$-minimal model
of $\mgn{\infty}$'', continuing the abuse of notation adopted for
\coh. Since a minimal model of a space may be obtained as a union
of $k$-minimal models, we call this union a ``minimal model of
$\mgn{\infty}$''. Each of these $k$-minimal models is $k$-formal in a
natural sense, see below, and the formality of $\mgn{\infty}$ follows.

First of all, we recall basic notions on $k$-minimal models, see
\cite{griffiths-morgan,lehmann,morgan}. From now on we will assume that our
spaces and algebras are connected and simply connected, {i.e.}, their
$H^0 = \nc$ and $H^1 = 0$. In application to moduli spaces, this is
the case as long as $g \ge 1$, see Harer \cite{harer2}. A DG
algebra $\gtm$ is called \emph{minimal} if it is free as a DG
commutative algebra, $\gtm^1 = 0$, and $d(\gtm) \subset \gtm^+ \cdot
\gtm^+$, where $\gtm^+ = \bigoplus_{i >0} \gtm^i$. A \emph{minimal
model} of a DG algebra $\A$ is a minimal DG algebra $\gtm$ along with
a
\emph{quasi-isomorphism} $\gtm \to \A$, a morphism of DG algebras
inducing an isomorphism on \coh. Every DG algebra $\A$ has a minimal
model, unique up to an isomorphism, which is in its turn unique up to
homotopy.

Let $k \ge 0$ be an integer. A $k$-\emph{minimal model} of a DG
algebra $\A$ is a minimal algebra $\gtm(k)$ generated by elements in
degrees $\le k$ along with a morphism $\gtm(k) \to \A$ inducing an
isomorphism on \coh\ in degrees $\le k$ and an injection in degree
$k+1$. A $k$-minimal model is unique up to an isomorphism uniquely
defined up to homotopy. If one has an increasing sequence of
embeddings
\[
\gtm(0) \subset \gtm(1) \subset \gtm(2) \subset \dots
\]
together with morphisms $\gtm(k) \to \A$ compatible with each other,
so that $\gtm(k)$ is a $k$-minimal model of $\A$, then the union $\gtm
=
\bigcup_k \gtm(k)$ along with the natural morphism $\gtm \to \A$ is a
minimal model of $\A$.

We will call a DG algebra $\A$ (or a space whose rational homotopy
type we consider) \emph{formal} if a minimal model $\gtm$ of $\A$ is
isomorphic to a minimal model of its \coh\ $H^\bullet(\A)$ taken with
the zero differential. This is equivalent to saying that there is a
quasi-isomorphism $\gtm \to H^\bullet(\A)$. In this case the rational
homotopy type of $\A$ (or the space) is determined by its rational
\coh\ ring. Formality implies that all Massey products are zero. If
the \coh\ algebra $H^\bullet(\A)$ is free as a graded commutative
algebra, then $\A$ is formal; see Proposition~\ref{k-formal} for a
$k$-version of this statement.  Another example of a formal space is
any compact K\"ahler manifold, a famous result of Griffiths, Deligne,
Morgan, and Sullivan \cite{dgms}.  We will similarly call a DG algebra
$\A$ or a space $k$-\emph{formal}, if $k$-minimal models of $\A$ and
its \coh\ are isomorphic.

We also need to prove the following proposition establishing a
particular case of $k$-formality. We say that a graded commutative
algebra $\C$ is \emph{$k$-free} if there exists a graded vector space
$V = \bigoplus_{i=0}^k V^k$ and a mapping $V \to \C$ of graded vector
spaces defining a morphism $S(V) \to \C$ of graded algebras, where
$S(V)$ is the free DG commutative algebra on $V$, which is an
isomorphism in degrees $\le k$ and an injection in degree $k+1$.

\begin{prop}
\label{k-formal}
Suppose that the \coh\ $H^\bullet(\A)$ of a DG algebra $\A$ is
$k$-free, based on a graded vector space $V$. Then $S(V)$ is a
$k$-minimal model of $\A$, and $\A$ is $k$-formal.
\end{prop}

\begin{proof}
If $H^\bullet(\A)$ is $k$-free, then there exists a linear injection
$\phi: V \hookrightarrow H^\bullet(\A)$. Pick a linear mapping $V \to
\A$ which takes each element $v$ of $V$ to a cocycle representing the
\coh\ class $\phi(v)$. Then the natural morphism of graded commutative
algebras $S(V) \to \A$ obviously respects the differentials and
satisfies the axioms of a $k$-minimal model of $\A$. By assumption,
$S(V)$ is at the same time a $k$-minimal model of $H^\bullet(\A)$,
whence $k$-formality of $\A$.
\end{proof}

Now we are ready to present the main result of the paper.

\begin{thm}
\begin{enumerate}
\item The moduli space $\mgn{g}$ is $k$-formal for $g \ge 2k+3$.
\item The subalgebra $H^\bullet (\mgn{\infty}) (k)$ of
$H^\bullet (\mgn{\infty})$ generated in degrees $\le k$ is a
$k$-minimal model of $\mgn{g}$ for $g \ge 2k+3$.
\item A $k$-minimal model of $\mgn{g}$ is independent of $g$ as long as
$g \ge 2k+3$. We will call it a $k$-\emph{minimal model of}
$\mgn{\infty}$.
\item The $k$-minimal models of $\mgn{\infty}$ form am increasing sequence
of embeddings. The union, a \emph{minimal model of} $\mgn{\infty}$, is
isomorphic to its \coh\ $H^\bullet(\mgn{\infty})$. In particular,
$\mgn{\infty}$ is formal.
\end{enumerate}
\end{thm}

\begin{proof}
1. The stable \coh\ $H^\bullet(\mgn{\infty})$ is a free graded
commutative algebra, \emph{i.e.}, isomorphic to $S(V)$ for a graded
vector space $V$, according to Miller-Morita's theorem \cite{miller}
for $n=0$ and Looijenga's handling \cite{looijenga} of the $n \ge 0$
case. Moreover, $H^0(\mgn{\infty}) = \nc$ and $H^1(\mgn{\infty}) = 0$
for $g \ge 1$, see Harer \cite{harer2}. Given a nonnegative integer
$k$, the groups $H^k(\mgn{g})$ are known to stabilize as soon as $g
\ge 2k+1$; this is the Harer-Ivanov Stability Theorem
\cite{harer:stab,ivanov}. Therefore, $H^\bullet(\mgn{g})$ is $k$-free
for $g \ge 2k+3$: the mapping $V^{\le k} \to H^\bullet (\mgn{g})$
makes $H^\bullet (\mgn{g})$ a $k$-free graded algebra.
Proposition~\ref{k-formal} then implies that $\mgn{g}$ is $k$-formal
for $g \ge 2k+3$, $S(V^{\le k})$ being a $k$-minimal model of it.

2 and 3. Notice that since $H^\bullet(\mgn{\infty})$ is free, the
subalgebra $H^\bullet(\mgn{\infty})(k)$ generated in degrees $\le k$
can be identified with $S(V^{\le k})$, which we have just seen to be a
$k$-minimal model of $\mgn{g}$ for $g \ge 2k+3$.

4. The subalgebras
\[
H^\bullet (\mgn{\infty}) (0) \subset H^\bullet (\mgn{\infty}) (1)
\subset H^\bullet (\mgn{\infty}) (2) \subset \dots
\]
form an increasing sequence of $k$-minimal models of $\mgn{\infty}$,
therefore, their union, $H^\bullet (\mgn{\infty})$, is a minimal model
of $\mgn{\infty}$, and thereby $\mgn{\infty}$ is formal.
\end{proof}

In view of this result, Mumford's Conjecture \cite{mumford}, if true,
implies the following refinement:
\emph{the polynomial algebra on the Mumford-Morita-Miller classes
$\kappa_i$, $i=1, 2, \dots,$ and the first Chern classes $c_1(T_i)$ of
the ``tangent at the $i$th puncture'' bundles, $i = 1, 2 , \dots, n$,
with a zero differential is the stable minimal model of the moduli
space $\mgn{g}$ as $g \to \infty$}.

\bibliographystyle{alpha}

\begin{thebibliography}{DGMS75}

\bibitem[Del71]{del:hodge}
P.~Deligne.
\newblock Th\'eorie de {H}odge. {I}{I}.
\newblock {\em Inst. Hautes \'Etudes Sci. Publ. Math.}, (40):5--57, 1971.

\bibitem[DGMS75]{dgms}
P.~Deligne, P.~A. Griffiths, J.~W. Morgan, and D.~Sullivan.
\newblock Real homotopy theory of {K}\"ahler manifolds.
\newblock {\em Invent. Math.}, 29(3):245--274, 1975.

\bibitem[GJ94]{gj}
E.~Getzler and J.~D.~S. Jones.
\newblock Operads, homotopy algebra and iterated integrals for double loop
  spaces.
\newblock Preprint, Department of Mathematics, MIT, March 1994.
\newblock \texttt{hep-th/9403055}.

\bibitem[GM81]{griffiths-morgan}
P.~A. Griffiths and J.~W. Morgan.
\newblock {\em Rational homotopy theory and differential forms}, volume~16 of
  {\em Progress in Mathematics}.
\newblock Birkh\"auser, Boston, MA, 1981.

\bibitem[Har85]{harer:stab}
J.~L. Harer.
\newblock Stability of the homology of the mapping class groups of orientable
  surfaces.
\newblock {\em Ann. of Math. (2)}, 121(2):215--249, 1985.

\bibitem[Har88]{harer2}
J.~L. Harer.
\newblock The \coh\ of the moduli space of curves.
\newblock In {\em Theory of moduli (Montecatini Terme, 1985)}, volume 1337 of
  {\em Lecture Notes in Math.}, pages 138--221. Springer, Berlin, 1988.

\bibitem[HS79]{halperin-stasheff}
S.~Halperin and J.~Stasheff.
\newblock Obstructions to homotopy equivalences.
\newblock {\em Adv. in Math.}, 32(3):233--279, 1979.

\bibitem[Iva93]{ivanov}
Nikolai~V. Ivanov.
\newblock On the homology stability for {T}eichm\"uller modular groups: closed
  surfaces and twisted coefficients.
\newblock In {\em Mapping class groups and moduli spaces of Riemann surfaces
  (G\"ottingen, 1991/Seattle, WA, 1991)}, volume 150 of {\em Contemp. Math.},
  pages 149--194. Amer. Math. Soc., Providence, RI, 1993.

\bibitem[Kon94]{maxim}
M.~Kontsevich.
\newblock Feynman diagrams and low-dimensional topology.
\newblock In {\em First European Congress of Mathematics, Vol.\ II (Paris,
  1992)}, volume 120 of {\em Progr. Math.}, pages 97--121. Birkh\"auser, Basel,
  1994.

\bibitem[KSV95]{ksv1}
T.~Kimura, J.~Stasheff, and A.~A. Voronov.
\newblock On operad structures of moduli spaces and string theory.
\newblock {\em Commun. Math. Phys.}, 171:1--25, 1995.
\newblock \texttt{hep-th/9307114}.

\bibitem[KSV96]{ksv2}
T.~Kimura, J.~Stasheff, and A.~A. Voronov.
\newblock Homology of moduli of curves and commutative homotopy algebras.
\newblock In {\em The Gelfand Mathematical Seminars, 1993--1995}, Gelfand Math.
  Sem., pages 151--170. Birkh\"auser Boston, Boston, MA, 1996.
\newblock \texttt{alg-geom/9502006}.

\bibitem[KVZ97]{kvz}
T.~Kimura, A.~A. Voronov, and G.~J. Zuckerman.
\newblock Homotopy {G}erstenhaber algebras and topological field theory.
\newblock In {\em Operads: Proceedings of Renaissance Conferences (Hartford,
  CT/Luminy, 1995)}, volume 202 of {\em Contemp. Math.}, pages 305--333,
  Providence, RI, 1997. Amer. Math. Soc.
\newblock \texttt{q-alg/9602009}.

\bibitem[Leh90]{lehmann}
D.~Lehmann.
\newblock Th\'eorie homotopique des formes diff\'erentielles (d'apr\'es {D}.
  {S}ullivan).
\newblock {\em Ast\'erisque}, (45):145, 1990.

\bibitem[Loo96]{looijenga}
E.~Looijenga.
\newblock Stable \coh\ of the mapping class group with symplectic
  coefficients and of the universal {A}bel-{J}acobi map.
\newblock {\em J. Algebraic Geom.}, 5(1):135--150, 1996.

\bibitem[Mil86]{miller}
E.~Y. Miller.
\newblock The homology of the mapping class group.
\newblock {\em J. Differential Geom.}, 24(1):1--14, 1986.

\bibitem[Mor78]{morgan}
J.~W. Morgan.
\newblock The algebraic topology of smooth algebraic varieties.
\newblock {\em Inst. Hautes \'Etudes Sci. Publ. Math.}, (48):137--204, 1978.

\bibitem[Mum83]{mumford}
D.~Mumford.
\newblock Towards an enumerative geometry of the moduli space of curves.
\newblock In {\em Arithmetic and geometry, Vol. II}, volume~36 of {\em Progr.
  Math.}, pages 271--328. Birkh\"auser Boston, Boston, Mass., 1983.

\bibitem[Pik95]{pikaart}
M.~Pikaart.
\newblock An orbifold partition of $\overline {M}{}\sp n\sb g$.
\newblock In {\em The moduli space of curves (Texel Island, 1994)}, volume 129
  of {\em Progr. Math.}, pages 467--482. Birkh\"auser Boston, Boston, MA, 1995.

\bibitem[Sta63]{jim}
J.~Stasheff.
\newblock Homotopy associativity of {$H$}-spaces. {I}, {I}{I}.
\newblock {\em Trans. Amer. Math. Soc.}, 108:275--312, 1963.

\bibitem[Sul77]{sullivan}
D.~Sullivan.
\newblock Infinitesimal computations in topology.
\newblock {\em Inst. Hautes \'Etudes Sci. Publ. Math.}, (47):269--331, 1977.

\end{thebibliography}


\end{document}